\newif\ifsmf
\IfFileExists{smfart.cls}{\smftrue}{\smffalse}

\ifsmf
  \documentclass{smfart}
\else
  \documentclass[11pt,a4paper]{amsart}
\fi

\usepackage[T1]{fontenc}
\usepackage{mathptm}
\usepackage[francais]{babel}
\usepackage{mathrsfs}

\usepackage[matrix,arrow,curve,cmtip]{xy}
\NoCompileMatrices

\theoremstyle{plain}
  \newtheorem{theo}[subsection]{Th\'eor\`eme}
  \newtheorem{prop}[subsection]{Proposition}
  \newtheorem{lemm}[subsection]{Lemme}
  \newtheorem{coro}[subsection]{Corollaire}
\theoremstyle{remark}
  \newtheorem{defi}[subsection]{Definition}
  \newtheorem{rema}[subsection]{Remarque}
  \numberwithin{equation}{subsection}
  \def\theequation{\arabic{section}.\arabic{subsection}.\arabic{equation}}

\def\Acris{A_\cris}
\def\Ainf{A_{\text{\upshape inf}}}

\def\BdR{\mathbf{B}_{\text{\upshape dR}}}
\def\Bst{\mathbf{B}_{\text{\upshape st}}}
\def\C {\mathbf C}
\def\D {\mathbf D}
\def\Q {\mathbf Q}
\def\Z {\mathbf Z}
\def\Cris{{\text{\upshape Cris}}}
\def\CalExt{\operatorname{\text{$\mathcal E$\kern-0.2em\itshape xt}}}
\def\CalHom{\operatorname{\text{$\mathcal H$\kern-0.2em\itshape om}}}
\def\RHom{\mathop{\mathbf{R}\CalHom}}
\def\cris{{\text{\upshape cris}}}
\def\Fil{\operatorname{Fil}}
\def\ga{\mathbf{G}_{\mathrm{a}}}
\def\gm{\mathbf{G}_{\mathrm{m}}}
\def\PN{{\og $p^n$\fg}}
\let\hra\hookrightarrow
\let\ra\rightarrow
\let\lra\longrightarrow
\def\Ext {\operatorname{Ext}}
\def\Gal {\operatorname{Gal}}
\def\Hom {\operatorname{Hom}}
\def\Spec {\operatorname{Spec}}
\def\ker {\operatorname{Ker}}
\let\eps\varepsilon
\let\phi\varphi
\let\mathcal\mathscr

\begin{document}

\ifsmf

\alttitle{Crystalline Dieudonn\'e theory and $p$-adic periods}
\begin{altabstract}
  We propose in this paper a theory of $p$-adic periods for finite
  flat group schemes. To this aim, we use the crystalline Dieudonn\'e
  theory of Berthelot, Breen and Messing, together with the crystalline
  interpretation of Fontaine's rings.
\end{altabstract}

\fi

\title[Th\'eorie de Dieudonn\'e cristalline et p\'eriodes $p$-adiques]
      {Th\'eorie de Dieudonn\'e cristalline \\ et p\'eriodes $p$-adiques}

\author{Antoine Chambert-Loir}
\address{Institut de Math\'ematiques de Jussieu \\
         Universit\'e Pierre et Marie Curie \\
         4 place Jussieu, F-75252 Paris Cedex 05}
\email{chambert@math.jussieu.fr}
\urladdr{http://www.math.jussieu.fr/~chambert}
\date{1\raisebox {.5ex}{\tiny er} septembre 1997,
      version corrig\'ee, 22 juin 1998, puis 13 ao\^ut 1998 }
\thanks{Paru dans \emph{Bull. Soc. Math. Fr.} \textbf{126} (1998), no. 4, 545-562}
\keywords{Th\'eorie de Dieudonn\'e, p\'eriodes $p$-adiques,
  cohomologie cristalline}
\subjclass{14F30, 14L}

\begin{abstract}
Nous proposons dans ce texte une th\'eorie des p\'eriodes $p$-adiques pour
des sch\'emas en groupes finis localement libres. Nous utilisons pour ce faire
la th\'eorie de Dieudonn\'e cristalline de Berthelot, Breen et Messing,
ainsi que l'interpr\'etation cristalline des anneaux de Fontaine.
\end{abstract}

\maketitle

\setcounter{tocdepth}{1}
\tableofcontents

\section{Introduction} \label{sec:intro}

Le but de cet article est de montrer qu'il existe pour les sch\'emas en
groupes finis localement libres
un th\'eor\`eme de comparaison $p$-adique {\og cristallin\fg}
analogue \`a celui que Fontaine a mis en \'evidence pour la
d\'ecomposition de Hodge--Tate.

Pour fixer les notations, soient $p$ un nombre premier,
$K$ un corps $p$-adique,
$\mathfrak O_K$ son anneau d'entiers~; on se donne une
$\mathfrak O_K$-alg\`ebre $R_0$ de caract\'eristique~$0$,
int\`egre et $p$-adiquement compl\`ete
et $R$ d\'esigne le compl\'et\'e $p$-adique de la cl\^oture int\'egrale de
$R_0$ dans une cl\^oture alg\'ebrique du corps des fractions de $R_0$.
Alors, $\Omega^1_{R/R_0}$ d\'esigne le module des diff\'erentielles
de K\"ahler continues de $R$ sur $R_0$.
Enfin, soit $A$ un sch\'ema ab\'elien sur $R_0$~; on
note $A'$ son sch\'ema ab\'elien dual, ainsi que $\omega_A$ le $R_0$-module
des formes diff\'erentielles invariantes sur $A$.

On peut alors r\'esumer ainsi une partie
des accouplements de p\'eriodes
pour les vari\'et\'es ab\'eliennes (ou les groupes de
Barsotti--Tate, voire les groupes formels)~:
\def\mylab#1{{\itshape #1. ---}\ \ignorespaces}
\begin{enumerate}
\def\theenumi{\alph{enumi}}
\item
\mylab{Formes diff\'erentielles de premi\`ere esp\`ece et diff\'erentielles
de K\"ahler}
Il s'agit d'une application $\omega_A\times T_p(A(R))\lra
T_p(\Omega^1_{R/R_0})$ mise en \'evidence par J.-M.~Fontaine.
Elle est explicit\'ee dans~\cite{fontaine82} quand $R_0=\mathfrak O_K$~;
dans ce cas, $T_p(\Omega^1_{R/R_0})$ s'identifie \`a $\C_p(1)$
comme $\Gal(\bar K/K)$-module.
\item
\mylab {Formes de seconde esp\`ece modulo les formes de premi\`ere esp\`ece}
Introduite par R.~Coleman~\cite{coleman84},
c'est une application $T_p(A(R))\lra \omega_{A'}\otimes_{R_0}R$ 
construite \`a l'aide de 
l'extension universelle de $A$ par un groupe vectoriel
(cf.~\cite{mazur-m76}).
\item
\mylab{Cohomologie cristalline ou de de Rham}
Plusieurs auteurs 
(citons notamment Fontaine--Messing,
Colmez~\cite{colmez92},
Candilera--Cristante~\cite{candilera-c95},
Wintenberger~\cite{wintenberger94})
cons\-truisent
diverses applications bilin\'eaires de la forme
$ D(A)\times T_p(A(R))\lra B$, o\`u $D(A)$ est un $R_0$-module
canoniquement attach\'e \`a $A$ (typiquement,
cohomologie cristalline ou de de Rham), et $B$ une $R_0$-alg\`ebre
du type de celles introduites par Fontaine~\cite{fontaine82,fontaine94}.
\end{enumerate}
De plus, les deux premiers accouplements expos\'es proviennent
d'une th\'eorie analogue pour les {\em sch\'emas en groupes finis 
localement libres}.
Si $G$ est un $R_0$-sch\'ema en groupes fini localement libre, on a
alors des accouplements~:
\begin{enumerate}
\def\theenumi{\alph{enumi}$'$}
\item
\mylab{Diff\'erentielles de K\"ahler et formes diff\'erentielles
invariantes}
Un $R$-point $x:\Spec R\ra G$ fournit par \'evaluation une
application $x^*:\omega_G\ra \Omega^1_{R/R_0}$, d'o\`u une application
bilin\'eaire $G(R)\times \omega_G\lra \Omega^1_{R/R_0}$
(Fontaine, cf.~\cite{fontaine82}).
\item
\mylab{Application universelle d'un sch\'ema en groupes fini localement
libre
dans un faisceau quasi-coh\'erent}
Mazur--Messing~\cite{mazur-m76} ont prouv\'e qu'il existe une application
universelle de $G$ dans un faisceau quasi-coh\'erent, de la forme
$\alpha_G:G\lra \omega_{G'}$, o\`u $G'$ est le dual de Cartier de $G$,
le lien avec l'accouplement b) a \'et\'e remarqu\'e pour la premi\`ere
fois par Crew~\cite{crew90}.
\end{enumerate}

Le but de ce texte est de montrer qu'au moins dans
certains cas, la th\'eorie~(c) correspond \`a un accouplement
analogue~(c$'$) pour les sch\'emas en groupes finis.
L'anneau de p\'eriodes sera,
si le groupe est tu\'e par $p^n$, une variante
de l'anneau $\Acris/p^n$ 
Quant au module de Dieudonn\'e, comme l'indique le
titre de l'expos\'e, il est fourni par la th\'eorie de Dieudonn\'e
cristalline de Berthelot, Breen et Messing. En fait, la construction de
l'accouplement provient de  constructions internes \`a la th\'eorie de
Dieudonn\'e cristalline~; calcul\'ees sur l'\'epaississement universel \`a
puissances divis\'ees, on obtient l'accouplement de p\'eriodes.

Faltings d\'emontre dans~\cite{faltings90} un th\'eor\`eme de comparaison
entre sch\'emas en groupes finis localement libres et certains modules filtr\'es
munis de connexions convenables. S'il est tr\`es vraisemblablement
plus g\'en\'eral, le lien pr\'ecis entre ce r\'esultat et ceux de
cet article ne m'appara\^\i t pas clairement.

\medskip
{\em \footnotesize
Ces r\'esultats ont \'et\'e expos\'es au Colloque
{\og Probl\`emes de coefficients en cohomologie cristalline et en cohomologie
rigide\fg} qui s'est tenu les 28--30 avril 1997 \`a l'Institut Henri
Poincar\'e, Paris. Je tiens \`a remercier Pierre Berthelot
pour l'int\'er\^et qu'il a port\'e \`a ce travail
en m'y invitant.
\par}

\section{Th\'eorie de Dieudonn\'e cristalline} \label{sec:tdc}

\subsection{}
Concernant la th\'eorie de Dieudonn\'e cristalline, nous utiliserons
les notations, m\'ethodes et r\'esultats
\'elabor\'es dans~\cite{berthelot-b-m82}.
Soit $S$ un sch\'ema plat sur $\Sigma=\Spec\Z_p$ (\'eventuellement un
sch\'ema formel $p$-adique sans $p$-torsion) et notons pour $n\geq 0$,
$S_n=V(p^n)\subset S$ et $\Sigma_n=\Spec \Z/p^n\Z$. On
travaille avec le (gros) site cristallin $\Cris(S_n,\Sigma)$ dont les
objets sont les pd-\'epaississements $(U,T,\gamma)$ ($p$ \'etant localement
nilpotent sur $T$, voire seulement topologiquement nilpotent) munis de
morphismes $U\ra S_n$, $T\ra \Sigma$ compatibles aux puissances
divis\'ees sur les id\'eaux $\mathcal J_T=\ker(\mathcal O_T\ra\mathcal
O_U)$ et $p\Z_p\subset\Z_p$. On remarque que $(S_n,S)$ muni de ses
puissances divis\'ees canoniques sur $p\mathcal O_S\subset\mathcal O_S$
d\'efinit un objet de $\Cris(S_n,\Sigma)$
Soit $G$ un sch\'ema en groupes fini localement libre sur $S$, tu\'e par une
puissance de $p$. Son image $\underline G$ par l'immersion
$i_{S_n/\Sigma}$ est le faisceau ab\'elien sur $\Cris(S_n,\Sigma)$ dont
les sections sur $(U,T,\gamma)$ sont donn\'ees par $G(U)$. Berthelot,
Breen et Messing d\'efinissent alors dans~\cite{berthelot-b-m82}
un cristal de Dieudonn\'e $\D(G)$ pour $G$~;
c'est le faisceau sur le site cristallin des extensions
locales de $\underline G$ par le faisceau structural. Le module de
Dieudonn\'e que nous utiliserons est alors d\'efini par~:

\begin{defi}
On pose $D(G)=\D(G)_{(S_n,S)}=
\CalExt^1_{S_n/\Sigma}(\underline G,\mathcal O_{S_n/\Sigma})_{(S_n,S)}$.
\end{defi}
C'est un $\mathcal O_S$-module de pr\'esentation finie, tu\'e par une
puissance de $p$ (cf.~\cite{berthelot-b-m82}, 3.1.3).
Si $S$ est le spectre de l'anneau des entiers d'un
corps local dont l'indice de ramification est inf\'erieur \`a $p-1$, il
s'identifie au module de Dieudonn\'e de la fibre sp\'eciale de $G$
tensoris\'e par $\mathcal O_S$. Sauf si l'endomorphisme de Frobenius de
$S_n$ se rel\`eve \`a $S$, il n'y a {\em pas} de Frobenius.

\begin{rema}
De Jong utilise dans \cite{dejong93} le faisceau des extensions
locales de $\underline G$ par l'id\'eal \`a puissances divis\'ees canonique.
Le module de Dieudonn\'e calcul\'e ainsi se compare alors au module de
Dieudonn\'e classique dans le cas o\`u $S_n$ est le spectre d'un corps
parfait. D'apr\`es la proposition 7.1 de~\cite{dejong93}, notre module de
Dieudonn\'e s'identifie \`a celui \'etudi\'e dans loc.\ cit.\ pour le tordu par
Frobenius de $G\times_S S_n$.
Si $G$ est le noyau de la multiplication par $p^n$
d'un sch\'ema ab\'elien $A$, c'est cependant notre normalisation qui
correspond \`a la r\'eduction modulo $p^n$ du premier groupe de cohomologie
de de Rham de $A$, ce qui explique son choix
(cf. aussi~\cite{colmez93}, haut de la page 640).
\end{rema}

Supposons que $G$ soit annul\'e par $p^n$. La construction \`a la base de
notre th\'eorie est l'application not\'ee \PN
dans~\cite{berthelot-b-m82}, p.~173--174~:

\begin{defi}
L'application \PN associe \`a une extension
locale de $\underline G$ par $\mathcal O_{S_n/\Sigma}$ l'homomorphisme
local $\underline G\lra \mathcal O_{S_n/\Sigma}/p^n$ d\'efini par le
diagramme du serpent de la multiplication par $p^n$ dans l'extension.
\end{defi}
Autrement dit, on rel\`eve un point de $\underline G$ en un point de
l'extension puis on le multiplie par $p^n$, pour obtenir un \'el\'ement de
$\mathcal O_{S_n/\Sigma}$ bien d\'efini  modulo $p^n$. La similitude avec
la construction de Coleman rappel\'ee au a) de l'introduction est alors
claire. D'apr\`es le corollaire 4.2.9 de~\cite{berthelot-b-m82},
cette application,
\'evalu\'ee sur tout \'epaississement $(U,T)$ tel que $T$
est la r\'eduction modulo $p^n$ d'un sch\'ema formel
plat sur $\Z_p$,
est un isomorphisme.
(Cf. aussi la proposition 3.10 de \cite{dejong93},
o\`u il suffit que $T$ soit la
r\'eduction modulo $p^n$ d'un sch\'ema plat sur $\Z/p^{n+1}\Z$.)

\subsection{Dualit\'e}
La dualit\'e de Cartier des sch\'emas en groupes finis (ou
celle des sch\'emas ab\'eliens) s'incarne dans la th\'eorie de Dieudonn\'e
cristalline au moyen d'un morphisme 
\begin{equation}\label{eq:gm-o1}
\underline{\gm}\lra\mathcal O_{S_n/\Sigma}[1]
\end{equation}
(dans la cat\'egorie d\'eriv\'ee) d\'efini
dans~\cite{berthelot-b-m82}, 5.2.1, p. 213 et que nous rappelons ici.
Partons de la suite exacte canonique
de faisceaux sur $\Cris(S_n/\Sigma)$
$$ 0\lra 1+\mathcal J_{S_n/\Sigma} \lra \mathcal O_{S_n/\Sigma}^* \lra
\underline{\gm} \lra 0. $$
Comme l'id\'eal $\mathcal J_{S_n/\Sigma}$ est muni de puissances
divis\'ees, il existe une application logarithme $\log:1+\mathcal
J_{S_n/\Sigma}\lra\mathcal J_{S_n/\Sigma}$ d\'efini (\cite{berthelot-b-m82},
3.2.7.3, p.~129]) par la formule
$$ \log (1+x) = \sum_{n=1}^\infty (-1)^{n-1} (n-1)! x^{[n]} $$
(qui est une somme localement finie dans le cas localement nilpotent,
et une s\'erie convergente dans le cas topologiquement nilpotent).
Le \emph{push-out\/} de la suite exacte pr\'ec\'edente par le morphisme
$$ 1+\mathcal J_{S_n/\Sigma} \xrightarrow{\log} \mathcal
J_{S_n/\Sigma}\hookrightarrow \mathcal O_{S_n/\Sigma} $$
nous fournit une extension
\begin{equation} \label{eq:us}
0\lra \mathcal O_{S_n/\Sigma} \lra \mathcal U_{S_n/\Sigma}\lra
\underline{\gm} \lra 0, 
\end{equation}
d'o\`u le morphisme annonc\'e dans la cat\'egorie d\'eriv\'ee.

La dualit\'e de Cartier identifie $G'$ \`a $\Hom(G,\gm)$, et de m\^eme
apr\`es application de $i_{S_n/\Sigma}$, $\underline{G'}=\Hom(\underline
G,\underline{\gm})$ (les deux $\Hom$ \'etant consid\'er\'es dans la cat\'egorie
des faisceaux ab\'eliens sur les sites respectivement fppf et
cristallin). Son image par le morphisme~\eqref{eq:gm-o1}
nous fournit un homomorphisme de faisceaux ab\'eliens sur le site cristallin
\begin{equation} \label{eq:g'-dg}
\underline{G'} \longrightarrow
    \CalExt^1_{S_n/\Sigma} (\underline G,\mathcal O_{S_n/\Sigma})=\D(G). 
\end{equation}
\begin{lemm}
Ce morphisme induit un isomorphisme
\begin{equation} \label{eq:dual1}
\CalHom_{S_n/\Sigma}(\D(G),\mathcal O_{S_n/\Sigma}) \stackrel\sim\lra
\CalHom_{S_n/\Sigma}(\underline G', \mathcal O_{S_n/\Sigma}). 
\end{equation}
\end{lemm}
\begin{proof}
L'application~\eqref{eq:g'-dg} est obtenue en prenant le $H^0$
du morphisme en cat\'egorie d\'eriv\'ee
\begin{equation} \label{eq:g'-delta}
 \underline{G'} \longrightarrow \Delta(G)[1] 
\end{equation}
induit par~\eqref{eq:gm-o1}
apr\`es application de $t_{0]}\RHom(\underline G,\cdot)$.
(Rappelons que le complexe de Dieudonn\'e $\Delta(G)$ est d\'efini comme
$t_{1]}\RHom(\underline G,\mathcal O_{S_n/\Sigma})$.)
Ainsi, le morphisme~\eqref{eq:dual1} est obtenu en appliquant
les foncteurs
$H^0$ et $\CalHom(\cdot,\mathcal O_{S_n/\Sigma})$
au morphisme~\eqref{eq:g'-delta}.

D'apr\`es le th\'eor\`eme 5.2.7 de~\cite{berthelot-b-m82},
en appliquant $t_{1]}\RHom(\cdot,\mathcal O_{S_n/\Sigma})$
\`a ce morphisme, on obtient un isomorphisme
entre complexes de Dieudonn\'e~:
\begin{equation} \label{eq:dual-delta}
 \Delta(G)^\vee[-1] \stackrel\sim\lra \Delta(G'). 
\end{equation}
Apr\`es application de $H^0$, il en r\'esulte donc un isomorphisme
\begin{equation}
\CalHom(\D(G),\mathcal O_{S_n/\Sigma}) \stackrel\sim\lra
\CalHom(\underline{G'},\mathcal O_{S_n/\Sigma})
\end{equation}
dont nous allons prouver qu'il s'identifie au morphisme~\eqref{eq:dual1}.
Autrement dit, il faut d\'emontrer que l'on peut sans
dommage \'echanger l'ordre des deux op\'erations
$(\mathbf R)\CalHom(\cdot,\mathcal O_{S_n/\Sigma})$
et $H^0$ \`a effectuer.

Or, consid\'erons la cat\'egorie $\mathcal C_0$ des complexes
$K^\cdot= [\cdots\ra K_1\ra K_0\ra 0\ra\cdots]$
dans $D(\mathcal O_{S_n/\Sigma})$
dont les composantes sont nulles en degr\'es~$>0$ et acycliques
pour $\Hom(\cdot,\mathcal O_{S_n/\Sigma})$~;
on a dans ce cas un isomorphisme canonique 
$$ H^0\RHom(K^\cdot,\mathcal O_{S_n/\Sigma}) \stackrel\sim\lra
\CalHom(H^0(K^\cdot),\mathcal O_{S_n/\Sigma}). $$
Notant $^\vee$ pour
$\CalHom(\cdot,\mathcal O_{S_n/\Sigma})$, 
les deux membres s'identifient en effet
tous deux \`a $\ker(K_0^\vee\ra K_{-1}^\vee)$.
Par suite, cet isomorphisme est fonctoriel pour les morphismes
dans la cat\'egorie d\'eriv\'ee qui se r\'ealisent comme
de morphismes de la cat\'egorie $\mathcal C_0$.

Or, $\Delta(G)$ admet, au moins localement sur $S_n$, 
une r\'esolution de longueur~$2$ par des $\mathcal O_{S_n/\Sigma}$-modules
libres de rang fini, concentr\'ee en rangs~$0$ et $1$, si bien
que $\Delta(G)[1]$ est bien un objet de $\mathcal C_0$.
De m\^eme, $G'$ admet des r\'esolutions par des (produits de) groupes
ab\'eliens libres du type $\Z[G^{\prime i}]$, ce qui entra\^{\i}ne
de plus que le morphisme~\eqref{eq:g'-delta} se r\'ealise comme un morphisme
de $\mathcal C_0$.
\end{proof}

\begin{coro} \label{coro:dualite}
Sur un objet $(U,T)$ du site cristallin o\`u l'homomorphisme \PN est
un isomorphisme, on dispose d'un isomorphisme de dualit\'e
$$ \D(G)_{(U,T)}^\vee \stackrel\sim\lra \CalHom_{S_n/\Sigma}
(\underline{G'}, \mathcal O_{S_n/\Sigma}) _{(U,T)}
\xleftarrow{\text{\PN}} \D(G')_{(U,T)}. $$
\end{coro}
En particulier, on retrouve
l'isomorphisme
$D(G)^\vee \simeq D(G')$ de~\cite{dejong93}, Proposition 6.3, 
dans lequel la dualit\'e $^\vee$ est d\'efinie par
$$D(G)^\vee=\CalHom_{\mathcal O_S}(D(G),\mathcal O_S/p^n). $$
Il est maintenant relativement ais\'e de relier l'homomorphisme \PN
aux dualit\'es sur les sch\'emas en groupes et sur les modules de
Dieudonn\'e. On a en effet le th\'eor\`eme~:

\begin{theo} \label{theo:dualite}
Supposons que $G$ est annul\'e par $p^n$ et
soit $(U,T)$ un objet du site cristallin
tel que l'homomorphisme \PN soit un isomorphisme.
Le diagramme suivant est commutatif
$$ \xymatrix{
{G_U\times G'_U } \ar[r]^-{\text{\PN}_{(U,T)}} 
                  \ar[d] &
{(D(G)^\vee\times D(G')^\vee) \otimes \mathcal O_T}
                \ar[r]^-{\text{dualit\'e}} &
{\mathcal O_T/p^n} \ar[d] \\
{\mu_{p^n\,(U,T)}}  \ar[rr]^-{\text{\PN}_{(U,T)}} & &
{\mathcal O_T/p^n.}}
$$
dans lequel la fl\`eche horizontale est fournie par le morphisme \PN
dans l'extension $\mathcal U_{S_n/\Sigma}$ de~\eqref{eq:us}.
\end{theo}
\begin{proof}
Via la d\'efinition~\ref{coro:dualite} de l'homomorphisme de dualit\'e, la ligne
du haut du diagramme associe \`a un couple $(x,y)\in G\times G'$ la
multiplication par $p^n$ d'un rel\`evement de $x$ dans l'extension
$y^*\mathcal U_{S_n/\Sigma}$ de $\underline G$ par $\mathcal
O_{S_n/\Sigma}$. Par fonctorialit\'e, $[x,y]\in\mu_{p^n}$ d\'esignant l'image
de $(x,y)$ par l'homomorphisme de dualit\'e de Cartier, c'est aussi la
multiplication par $p^n$ d'un relev\'e de $[x,y]$ dans l'extension
$\mathcal U_{S_n/\Sigma}$, ce qui \'etablit le th\'eor\`eme.
\end{proof}

\section{Construction d'un accouplement de p\'eriodes} \label{sec:acc}
\subsection{}
Soient $K$ un corps $p$-adique, $\mathfrak O_K$ son anneau des entiers,
Soit $R$ une $\mathfrak O_K$-alg\`ebre int\`egre, sans $p$-torsion
et compl\`ete pour la topologie $p$-adique~;
notons $F$ son corps des fractions.

Soit $\tilde F_p$ une extension maximale de $F$ telle
que si $\tilde R_p$ est la cl\^oture int\'egrale de $R$
dans $\tilde F_p$, alors $\tilde R_p[1/p]$ est \'etale sur $R[1/p]$.
Soit $\tilde F$ une extension alg\'ebrique de $F$ contenant $\tilde F_p$ ;
on note alors $\tilde R$
la cl\^oture int\'egrale de $R$ dans $\tilde F$, ainsi que
$\widehat{\tilde F}$ et $\widehat{\tilde R}$
les compl\'et\'es $p$-adique
de $\tilde F$ et $\tilde R$ respectivement.

\subsection{\'Epaississement universel}
On fait enfin l'hypoth\`ese suivante~: \par
{{\normalfont ($*$)\quad}\em
l'homomorphisme de Frobenius de $\tilde R/p\tilde R$ est surjectif.
\par\noindent}
Si $R$ est petit au sens de Faltings, 
(voir~\cite{faltings90}, \S2)
cette hypoth\`ese est v\'erifi\'ee pour $\tilde R=\tilde R_p$.
Comme dans~\cite{fontaine94} 2.2.1, \cite{wintenberger94}
et~\cite{tsuji97} A.1.5, cette hypoth\`ese permet d'assurer
l'existence d'un $R$-pd-\'epaississement
universel $p$-adique de $\tilde R/p\tilde R$ que
nous noterons $A_{\cris}(\tilde R/R)$.

Cet \'epaississement universel admet la description explicite
suivante~: soit $\mathcal R$ la limite projective des alg\`ebres
$\tilde R/p\tilde R$ relativement aux morphismes d'\'el\'evation
\`a la puissance
$p$\up{\`eme}, qui est une alg\`ebre parfaite de caract\'eristique~$p$.
L'anneau $W(\mathcal R)$ est alors muni d'un homomorphisme canonique
$W(\mathcal R)\ra \widehat{\tilde R}$, surjectif en vertu de
l'hypoth\`ese~($*$), d'o\`u un homomorphisme toujours surjectif
$$ \theta: R \otimes_{\Z_p} W(\mathcal R) \lra \widehat{\tilde R}. $$
L'anneau $A_{\cris}(\tilde R/R)$ est
alors le compl\'et\'e $p$-adique de l'enveloppe \`a
puissances divis\'ees de l'id\'eal $J_\cris(\tilde R/R)$,
noyau de $\theta$,
compatibles avec les puissances divis\'ees canoniques sur l'id\'eal $(p)$.

Lorsque $R=\mathfrak O_K$ et $\tilde R=\mathfrak O_{\bar K}$,
nous noterons $A_{\cris,K}$ l'anneau obtenu. 
Cette construction est fonctorielle en le couple $(R,\tilde R)$,
d'o\`u une action de $\Gal(\tilde F/F)$, ainsi qu'une structure
de $A_{\cris,K}$-alg\`ebre sur $A_{\cris}(\tilde R/R)$.

\subsection{}
Soit $S=\Spec R$ et $G$ un sch\'ema en groupes
fini localement libre sur $S$, annul\'e par une puissance de~$p$.
Remarquons que comme $G$ est fini et plat sur $R$ et $\tilde R$
est int\'egralement close dans~$\tilde F$,
il y a bijection entre les $\tilde F$-points
de $G_F$ et le sections $\Spec\tilde R\ra G$.

%

\subsection{}
On applique maintenant la th\'eorie du~\S\ref{sec:tdc}~;
l'\'epaississement \`a puissances
divis\'ees $A_{\cris}(\tilde R/R)$ d\'efinit notamment un objet
du site $\Cris(S_n/\Sigma)$.
Comme le cristal de Dieudonn\'e $\D(G)$ est un cristal, on a un
isomorphisme naturel $A_{\cris}(\tilde R/R)\otimes_{R} D(G)\simeq
\D(G)_{A_{\cris}(\tilde R/R)}$.
Jointe \`a l'homomorphisme naturel (r\'eduction modulo~$p^n$
des coordonn\'ees) $G(\tilde R)\lra G(\tilde R/p^n\tilde R)$,
la th\'eorie du \S\ref{sec:tdc}
nous fournit alors un accouplement
$$ \Phi_{G,n} : 
    G(\tilde R) \times D(G)
        \lra A_{\cris}(\tilde R/R)/p^n A_{\cris}(\tilde R/R). $$
Il est fonctoriel \`a la fois en $G$ et en le couple $(R,\tilde R)$,
et commute donc aux actions de $\Gal(\tilde F/F)$
sur $G(\tilde F)$ et $A_{\cris}(\tilde R/R)$.
Pour ne pas fixer un entier $n$ tel que $p^n $ annule~$G$,
on peut diviser cet
accouplement par $p^n$ et obtenir ainsi une application bilin\'eaire
$$ \Phi_G : G(\tilde R) \times D(G) 
       \lra A_{\cris}(\tilde R/R)\otimes_{\Z_p} (\Q_p/\Z_p).$$

\subsection{Cas du groupe multiplicatif}
Fixons une fois pour toutes un g\'en\'erateur
$(\eps_1,\ldots,\eps_n,\ldots)$ de $\Z_p(1)$.
Supposons dans ce paragraphe que $G=\mu_{p^n}$~; son module de Dieudonn\'e
$D(G)$ s'identifie \`a $\mathfrak O_S/p^n$, avec une base naturelle $dt/t$
fournie par la
restriction \`a $\underline G$ de l'extension $\mathcal U_{S_n/\Sigma}$
(cf.~\cite{dejong93}, p.~106).
Soit alors $x_n\in G(\tilde R)\subset\mathfrak O_{\bar K}\subset \tilde R$
tel que $[p^n]x_n=0$, la
p\'eriode de $dt/t$ contre $x_n$ est obtenue ainsi~:
choisissons
$\widetilde{x_n}\in A_{\cris,K}\subset A_{\cris}(\tilde R/R)$
relevant $x_n$, alors,
$\widetilde{x_n}^{p^n}\in 1+J_{\cris,K}\subset 1+J_\cris(\tilde R/R)$
et la p\'eriode vaut
$\log(\widetilde{x_n}^{p^n})$. Avec les notations de~\cite{fontaine94},
1.5.4, si $x_n=\eps_n$,
on a $\widetilde{\eps_n}^{p^n}=\nu(\eps)$ modulo $p^n A_{\cris,K}\subset
p^n A_\cris(\tilde R/R)$,
si bien que la p\'eriode est $t_{p,K}=2i\pi=\log\nu(\eps)$
(qui provient de $A_{\cris,K}$ mais consid\'er\'ee dans $A_\cris(\tilde
R/R)/p^n$).

Le th\'eor\`eme suivant est une cons\'equence facile de ce calcul 
et de la compatibilit\'e~\ref{theo:dualite}
de l'accouplement de p\'eriodes aux diverses dualit\'es.
\begin{theo}
L'application lin\'eaire 
$$\Phi_G: A_{\cris}(\tilde R/R)\otimes_\Z G(\tilde R)
     \longrightarrow \Hom(D(G),A_{\cris}(\tilde R/R)\otimes(\Q_p/\Z_p)) $$
a un noyau et un conoyau tu\'es par~$2i\pi=t_{p,K}\in A_{\cris,K}\subset
A_{\cris}(\tilde R/R)$.
\end{theo}
Commen\c{c}ons par le noyau.
Si $x\in A_{\cris}(\tilde R/R)\otimes_\Z G(\tilde R)$
s'envoie sur~$0$ par $\Phi_G$, le r\'esultat de dualit\'e 2.7
implique que pour tout
$x'\in A_{\cris}(\tilde R/R)\otimes G'(\tilde R)$,
$[x,x']\in A_{\cris}(\tilde R/R)\otimes \mu_{p^{\infty}}$ a pour image~$0$
par l'accouplement de p\'eriodes sur le groupe multiplicatif.
Autrement dit,
$t\Phi_{\gm}([x,x'])=0\in A_{\cris}(\tilde R/R)\otimes(\Q_p/\Z_p)$
pour tout $x'$,
d'o\`u $t_{p,K}x=0$
dans $A_{\cris}(\tilde R/R)\otimes_\Z G(\tilde R)$.

Pour le conoyau,
soit $\omega\in\Hom(D(G),A_{\cris}(\tilde R/R)\otimes(\Q_p/\Z_p))$.
Alors, l'application $A_{\cris}(\tilde R/R)\otimes G'(\tilde R)
\lra A_{\cris}(\tilde R/R)/p^n$ d\'efinie par $x'\mapsto
t\langle\omega,\Phi_{G'}(x')\rangle$ provient par dualit\'e d'un \'el\'ement
$x\in A_{\cris}(\tilde R/R)\otimes G(\tilde R)$ qui v\'erifie
$\Phi_G(x)=t_{p,K}\omega$.
Le th\'eor\`eme est ainsi d\'emontr\'e.

\subsection{Remarques}
Ce th\'eor\`eme et sa d\'emonstration sont bien s\^ur inspir\'es
du th\'eor\`eme dit de {\og presque d\'ecomposition de Hodge--Tate\fg}
pour les sch\'emas en groupes finis et plats tu\'es par
une puissance de $p$ sur $\mathfrak O_K$,
th\'eor\`eme d\^u \`a J.-M. Fontaine (Corollaire du th\'eor\`eme 3
de~\cite{fontaine82}).

D'autre part, le rapporteur me signale un article
de G.~Faltings~\cite{faltings94} o\`u est construit
un homomorphisme similaire pour les groupes de Barsotti--Tate
sur $\mathfrak O_K$ : il est injectif et son conoyau est
annul\'e par $t_{p,K}$ (th\'eor\`eme 7 de \emph{loc. cit.}).

Lorsque $R=\mathfrak O_K$, \'etant donn\'e que tout $R$-sch\'ema en groupes
fini et plat se plonge dans un groupe de Barsotti--Tate,
on peut retrouver notre r\'esultat \`a partir de celui de Faltings.

R\'eciproquement, appliqu\'e aux noyaux des multiplications par $p^n$
dans un groupe de Barsotti--Tate,
notre th\'eor\`eme permet de retrouver celui de Faltings et, en fait,
de l'\'etendre \`a une base plus g\'en\'erale que $\Spec\mathfrak O_K$.

\subsection{Que vaut $t_p$ ?}
Le th\'eor\`eme pr\'ec\'edent ne serait d'aucun int\'er\^et
si $t_{p,K}$ \'etait nul dans $A_{\cris}(\tilde R/R)$ et il importe
en fait de conna\^{\i}tre la valuation $p$-adique de son annulateur 
$\mathfrak a(\tilde R/R)$ dans $A_{\cris}(\tilde R/R)\otimes (\Q_p/\Z_p)$. 

Lorsque $R=\mathfrak O_K$, J.-M.~Fontaine a calcul\'e dans~\cite{fontaine82}
l'annulateur de l'image de $t_{p,K}$ modulo $\Fil^2$. Il obtient
un id\'eal de $\mathfrak O_{\bar K}$
dont tout g\'en\'erateur a une valuation $p$-adique \'egale \`a
$\frac1{p-1} + v_p(\mathfrak d_{K/\Q_p})$, o\`u $\mathfrak d_{K/\Q_p}$
d\'esigne la diff\'erente de $K/\Q_p$.

Si $\Spec R$ admet un point $\mathfrak O_K$-rationnel, il
r\'esulte de la fonctorialit\'e des anneaux $A_{\cris}$ 
un morphisme d'\'evaluation $A_{\cris}(\tilde R/R)\ra A_{\cris,K}$
qui est un inverse \`a gauche de l'inclusion $A_{\cris,K}\subset
A_{\cris}(\tilde R/R)$. Ainsi, l'annulateur de $t_{p,K}$ dans
$A_{\cris}(\tilde R/R)\otimes(\Q_p/\Z_p)$ ne sera pas plus gros
que dans $A_{\cris,K}\otimes\Q_p/\Z_p$.

Dans le cas o\`u $\tilde R[1/p]$ est \'etale sur $R[1/p]$,
les diff\'erentielles de K\"ahler fournissent (cf.~\cite{fontaine94},
1.4.6, ainsi que la remarque 1.4.8, voir aussi le
paragraphe~\ref{par:ordre1} plus bas) un \'epaississement
du premier ordre de $\tilde R$, et donc un \'epaississement
\`a puissances divis\'ees (toujours du premier ordre).
L'annulateur de $t_{p,K}$
dans $A_{\cris}(\tilde R/R)\otimes(\Q_p/\Z_p)$ a ainsi
une valuation au plus \'egale \`a $1/(p-1)$ plus la diff\'erente de
la cl\^oture alg\'ebrique de $K$ dans $F$.

\section{Comparaisons} \label{sec:compa}

Dans cette section, nous voulons montrer les relations entre
notre accouplement de p\'eriodes
et les th\'eories esquiss\'ees dans l'introduction.
Nous conservons les notations des paragraphes pr\'ec\'edents.

\subsection{\'Epaississements du premier ordre} 
\label{par:ordre1}
On rappelle que $\Ainf(\tilde R/R)= R \otimes W(\mathcal R)$
est le $R$-\'epaississement infinit\'esimal $p$-adique universel de
$\widehat{\tilde R}$, et que $A_{\cris}(\tilde R/R)$
est le compl\'et\'e $p$-adique
de l'enveloppe \`a puissances
divis\'ees de l'id\'eal noyau de $\theta$ compatibles aux puissances
divis\'ees canoniques sur l'id\'eal $(p)$. Si $\Fil$ est la filtration
induite sur $\Ainf$ (resp.~$A_{\cris,K})$  par les puissances du noyau
de $\theta$, la th\'eorie des puissances divis\'ees
(\cite{berthelot74b}, 3.3.4) permet d'affirmer que l'application naturelle
$$ \Ainf/ \Fil^2 \Ainf \longrightarrow \Acris/\Fil^2
A_{\cris} $$
est un isomorphisme si bien que modulo $\Fil^2$, $A_{\cris}$ s'identifie
\`a l'\'epaississement infinit\'esimal universel du premier ordre.

Supposons que $\tilde R[1/p]$ est \'etale sur $R[1/p]$. 
Les diff\'erentielles de K\"ahler fournissent alors un \'epaississement
infinit\'esimal $p$-adique du premier ordre
(voir~\cite{fontaine94}, 1.4.8, et plus bas lorsque $R=\mathfrak O_K$)~;
il est d'ailleurs universel, mais 
c'est l'existence de cet \'epaississement du premier ordre
qui nous importe, plus que son universalit\'e.

\subsection{Formes diff\'erentielles invariantes}
Soit $H$ un $S$-sch\'ema en groupes lisse.
Il est construit dans~\cite{berthelot-b-m82}, 3.2, un isomorphisme
fonctoriel
$$ \omega_{H/S_n} \ra
    \CalExt^1(i_{S_n/S_n*}(H),\mathcal J_{S_n/S_n})_{(S_n,S_n)}. $$
Compte tenu des homomorphismes naturels
\begin{multline*}
\CalExt^1( i_{S_n/S_n*}(H) ,\mathcal J_{S_n/S_n})_{S_n}
\lra \CalExt^1 ( i_{S_n/S_n*}(H),\mathcal J_{S_n/S_n})_{(S_n,S_n)} \\
\lra  \CalExt^1 ( i_{S_n/S_n*}(H),\mathcal O_{S_n/S_n})_{(S_n,S_n)} 
\end{multline*}
et $\omega_{H}\ra \omega_{H/S_n}$, on en d\'eduit un homomorphisme
fonctoriel
$ \omega_H \ra \mathbf D(H)_{S_n} $.
Ces homomorphismes fournissent un homomorphisme $\omega_G \ra\mathbf D(G)_{S_n}$
(en passant par la cat\'egorie d\'eriv\'ee) si l'on dispose d'une r\'esolution 
de $G$ par des sch\'emas en groupes lisses.

A. J. De Jong donne dans~\cite{dejong93} une description assez explicite
de cet homomorphisme qui va nous \^etre utile.
Sur le site $\Cris(S_n/S_n)$, on dispose d'un faisceau
$\pi_{S_n/S_n}^* H$ dont les sections sur un \'epaississement $(U,T)$
sont simplement $H(T)$, ce qui a un sens puisque $T$ est un $S_n$-sch\'ema.
Les homomorphismes $H(T)\ra H(U)$ induisent une
suite exacte
\begin{equation}
0 \ra H(\mathcal J_{S_n/S_n}) \ra \pi_{S_n/S_n}^* H \ra \underline H \ra 0.
\end{equation}
(On a not\'e $H(\mathcal J_{S_n/S_n})$ le faisceau sur $\Cris(S_n/S_n)$
dont les sections sur $(U,T)$ sont donn\'ees par $\ker (H(T)\ra H(U))$,
la surjectivit\'e vient de ce que $H$ est lisse.)
Le lemme de Poincar\'e cristallin permet de d\'efinir un homomorphisme
d'int\'egration des formes diff\'erentielles 
\begin{equation} \label{equa:poincare}
 \omega_{H/S_n}
  \ra \Hom (H(\mathcal J_{S_n/S_n}), \mathcal J_{S_n/S_n})_{S_n} , 
\end{equation}
d'o\`u en combinant les deux \'equations pr\'ec\'edentes, un homomorphisme
$$
\omega_{H/S_n} \ra \CalExt^1( \underline H,\mathcal J_{S_n/S_n})_{S_n}
$$
et en particulier un homomorphisme
\begin{equation}
\omega_{H/S_n} \ra \CalExt^1(\underline H,\mathcal J_{S_n/S_n})_{S_n}.
\end{equation}
On peut choisir $H$ de sorte \`a disposer d'une immersion ferm\'ee
$G\hra H$,  par exemple $H=\operatorname{Mor}_{\text{Sch}}(G^*,\gm)$ convient.
On en d\'eduit par fonctorialit\'e un homomorphisme
$$
\omega_{G/S_n} \ra \CalExt^1(\underline G,\mathcal J_{S_n/S_n})_{S_n}
 \ra \CalExt^1(\underline G,\mathcal O_{S_n/S_n})_{S_n}.
$$
puis, si $G$ est annul\'e par $p^n$, un homomorphisme
\begin{equation}
\omega_G \lra D(G)
\end{equation}
dont A. J. De Jong affirme (\emph{loc. cit.}, remarque 4.4)
qu'il est \emph{l'oppos\'e} de l'homomorphisme analogue
construit par Berthelot, Breen et Messing (\emph{loc. cit.}, 3.2.6).

Les deux paragraphes pr\'ec\'edents et notre accouplement $\Phi_G$
nous fournissent ainsi une application bilin\'eaire
$G(\tilde R)\times\omega_G \ra \Ainf/ (p^n \Ainf+\Fil^2 \Ainf )$
dont l'image est par construction contenue dans
$\Fil^1(\Ainf )/(p^n \Ainf +\Fil^2 \Ainf )$, d'o\`u
une application bilin\'eaire
\begin{equation} \label{eq:1espece}
\Phi^1_G : G(\tilde R) \times\omega_G \lra \Omega^1_{\tilde R/R}. 
\end{equation}
Or, il existe une telle application bilin\'eaire \`a la fois canonique
et explicite,
dont Fontaine (\cite{fontaine82}, 4.7) a montr\'e l'utilit\'e
dans le cas $R=\mathfrak O_K$.
Elle est donn\'ee par l'\'evaluation~: si $\omega\in\omega_G$ et
$x\in G(\tilde R)$, on peut voir $x$ comme une section $\Spec\tilde R\ra G$
au-dessus de $\Spec R$,
d'o\`u une forme diff\'erentielle $x^*\omega\in\omega_{\tilde R/R}$.

\begin{prop}
L'accouplement~\eqref{eq:1espece} est \'egal \`a
celui construit par Fontaine dans~{\em loc.\ cit.}
\end{prop}
Explicitons d'abord la suite exacte du milieu de la
page~68 de~\cite{fontaine94} qui fournit une description de
$A_{\cris,K}/(p^n A_{\cris,K}+\Fil^2)$~:
\begin{multline*}
0\lra \Omega^1_{\mathfrak O_{\bar K}/\mathfrak O_K}[p^n] \ra \\
\lra \mathfrak O'_{\bar K}/p^n = A_{\cris,K}/(p^n A_{\cris,K}+\Fil^2A_{\cris,K})
\lra \mathfrak O_{\bar K}/p^n\lra 0. 
\end{multline*}
Dans cette formule, $\mathfrak O'_{\bar K}$ est l'ensemble des \'el\'ements
$x \in \mathfrak O_{\bar K}$ tels que $dx=0$ dans $\Omega^1_{\mathfrak
O_{\bar K}/\mathfrak O_K}$. La derni\`ere fl\`eche est induite par
l'inclusion naturelle, sa surjectivit\'e r\'esulte de ce que
$\Omega^1_{\mathfrak O_{\bar K}/\mathfrak O_K}$ est $p$-divisible.
Le noyau de cette fl\`eche est form\'e des $x\in\mathfrak O'_{\bar K}/p^n$
tels que $x\in p^n\mathfrak O_{\bar K}$~; on \'ecrit alors $x=p^n u$,
$u$ \'etant d\'etermin\'e modulo $\mathfrak O'_{\bar K}$ et on lui
associe la diff\'erentielle de K\"ahler
$du\in \Omega^1_{\mathfrak O_{\bar K}/\mathfrak O_K}$ (qui est tu\'ee par
multiplication par $p^n$).

On fixe un entier $n$ tel que $p^n$ annule $G$ et on se place
sur $\Cris(S_n,S_n)$.
D'apr\`es la description rappel\'ee ci-dessus de l'homomorphisme
$\omega_G\ra D(G)$, les p\'eriodes de $x\in G(\tilde R)$ sont calcul\'ees
de la fa\c con suivante.
On plonge $G$ dans le sch\'ema en groupes lisse des morphismes
de sch\'emas de $G^*$ dans $\gm$. Si $\mathcal A_{G}$ est l'alg\`ebre
affine de $G$, $x$ correspond \`a un \'el\'ement inversible
$a$ de l'alg\`ebre
$\mathcal A_G^\vee\otimes_{\mathfrak O_K} \mathfrak O_{\bar K} $.
On le {\og rel\`eve\fg} en un \'el\'ement inversible
$\tilde a \in\mathcal A_G^\vee\otimes _{\mathfrak O_K} \mathfrak O'_{\bar K} $
tel que $\tilde a=a(1+p^n u)$ pour
$u\in\mathcal A_G^\vee\otimes \mathfrak O_{\bar K}$
tel que $d\tilde a=0$, $d$ d\'esignant la d\'erivation
$$\mathcal A_G^\vee\otimes\mathfrak O_{\bar K}
  \ra\mathcal A_G^\vee\otimes \Omega^1_{\mathfrak O_{\bar K}/\mathfrak O_K}. $$
On calcule alors $\tilde a^{p^n}$ :
$$ \tilde a^{p^n} = a^{p^n} (1+p^n u)^{p^n}
 = 1 + \sum_{k=1}^{p^n} {\textstyle\binom{p^n}k} p^{nk} u^k 
 = 1 + p^n \sum_{k=0}^{p^n-1}
         \frac{p^n}{k+1}{\textstyle\binom{p^n-1}k} p^{nk} u^{k+1}.
$$
Cet \'el\'ement rel\`eve $1$ et la diff\'erentielle de K\"ahler associ\'ee \`a
$\tilde a^{p^n}-1$ est ainsi
$$  \sum_{k=0}^{p^n-1} p^n {\textstyle\binom{p^n-1}k} p^{nk} u^k \, du 
= p^n (1+p^n u)^{p^n-1}\, du 
= (1+p^n u)^{p^n}  \frac{p^n \, du}{1+p^n u} = - \frac{da}{a}
$$
puisque 
$$0=\frac{d\tilde a}{\tilde a} = \frac{da}{a}+\frac{p^n \, du}{1+p^n u} $$
et que $da$ est annul\'e par $p^n$.

D'apr\`es~\cite{fontaine82}, \S\,4.5, L'\'el\'ement
$da/a$ de $\mathcal A_G^\vee \otimes_{\mathfrak O_K} \Omega^1_{\mathfrak
O_{\bar K}/\mathfrak O_K} $
appartient en fait \`a un sous-module canoniquement identifi\'e
\`a $\Hom(\omega_G,\Omega^1_{\mathfrak O_{\bar K}/\mathfrak O_K})$. 
Comme nous nous somme plac\'es sur un \'epaississement du premier ordre,
l'homomorphisme d'int\'egration des formes diff\'erentielles~\ref{equa:poincare}
fournie par le lemme de Poincar\'e cristallin
s'incarne en la dualit\'e entre formes diff\'erentielles
et espace tangent, si bien que notre accouplement $\Phi^1_G$
associe \`a $x\in G(\mathfrak O_{\bar K})$ ce m\^eme homomorphisme.

Finalement, l'homomorphisme en question est lui-m\^eme identifi\'e
dans la proposition~6
de~\cite{fontaine82} \`a l'homomorphisme d'\'evaluation des formes
diff\'erentielles, d'o\`u la proposition.

\begin{rema}
Nous n'avons donn\'e la d\'emonstration que dans le cas $R=\mathfrak O_K$,
car c'est dans ce cas que se place Fontaine dans~\cite{fontaine82},
et aussi dans ce cas que la litt\'erature explicite l'\'epaississement
du premier ordre \`a l'aide des diff\'erentielles de K\"ahler.
Les r\'esultats de~\cite{fontaine82} que nous avons utilis\'es
n\'ecessitent uniquement le fait que l'alg\`ebre affine $\mathcal A_G$
soit libre, et s'\'etendent au cas o\`u elle est localement libre.
Quant \`a l'interpr\'etation de l'\'epaississement du premier ordre
\`a l'aide de diff\'erentielles de K\"ahler, la remarque~1.4.8
de~\cite{fontaine94} montre qu'elle est valable
si $\tilde R[1/p]$ est \'etale sur $R[1/p]$.

Suivant cette m\^eme remarque,
il aurait \'et\'e possible de faisceautiser tout l'article en introduisant
le site syntomique de $\Spec R$. La construction du \S\ref{sec:tdc}
s'interpr\`ete alors comme un homomorphisme
$G\ra\Hom_{\text{SYN}}(D(G),\mathcal O_{n,\cris})$, dont noyau et
conoyau sont comme auparavant tu\'es par un \'el\'ement canonique~$t_p$,
section globale du faisceau $\mathcal O_{n,\cris}$.
La comparaison avec les diff\'erentielles de K\"ahler serait alors apparue
comme restriction de cette application au site syntomique-\'etale de
$\Spec R$ (cf.~\cite{fontaine94}, remarque 1.4.8).
\end{rema}

\subsection{Formes de seconde esp\`ece}
La comparaison est ici plus formelle. En consid\'erant le site
$\Cris(S_n,\Sigma)$, on associe \`a tout \'el\'ement de $\CalExt(\underline
G,\mathcal O_{S_n/\Sigma})$ une extension dans $\CalExt(\underline
G,\ga)=i_{S_n/\Sigma *}\Ext_{S_n}(G,\ga)$. Du point de vue des p\'eriodes,
cela consiste \`a appliquer le morphisme~$\theta$. En plongeant $G$ dans un
sch\'ema ab\'elien~$A$, $\Ext(G,\ga)$ est induit par l'extension universelle
et l'identification avec la th\'eorie de Coleman~\cite{coleman84},
{\em Note added in proof}, est apparente.

\subsection{Th\'eorie de Wintenberger}

J.-P.~Wintenberger construit dans~\cite{wintenberger94} un accouplement
de p\'eriodes pour les sch\'emas ab\'eliens.
Dans le cas consid\'er\'e dans le pr\'esent article,
la construction de Wintenberger est facile 
puisque nous nous sommes cantonn\'es
au cas de bonne r\'eduction. En effet, avec les notations de {\em loc.
cit.,} le point important est de montrer l'existence d'un sous-groupe
born\'e des $\Ainf[1/p]$-points de l'extension universelle,
le reste de la construction proc\'edant comme chez Coleman
par rel\`evement et multiplication par des puissances de $p$.
Or, dans le cas o\`u $A/R$ est un sch\'ema ab\'elien, le sous-groupe
des $\Ainf$-points de l'extension universelle convient visiblement,
ce qui montre que notre accouplement redonne celui de Wintenberger
par passage \`a la limite.

Au passage, cela prouve la non-d\'eg\'enerescence de l'accouplement de
Wintenberger, malheureusement non explicit\'ee dans~\cite{wintenberger94}.

\vskip -1cm
\section{Remarques et perspectives~?} \label{sec:remarks}

Tout d'abord, remarquons que les m\'ethodes d\'evelopp\'ees dans cet article
se retrouvent en de nombreux endroits de la litt\'erature sous une forme
souvent tr\`es proche, notamment dans les articles traitant de la plein
fid\'elit\'e du foncteur de Dieudonn\'e cristallin.

Compte tenu du travail~\cite{trihan96} de F.~Trihan concernant la th\'eorie
de Dieudonn\'e cristalline de niveau variable, il me semble aussi
que les m\'ethodes et r\'esultats expos\'es dans cet article s'\'etendent
au cas des cristaux de Dieudonn\'e de niveau fini.
Toutefois, les pd-\'epaississements universels
$p$-adiques de niveau quelconque n'ont \`a ma connaissance
pas b\'en\'efici\'e d'une \'etude aussi approfondie que  pour
le niveau~$0$.

D'autre part, et plus s\'erieusement,
nous ne prenons pas en compte la mauvaise r\'eduction,
ce qui est par exemple la principale source de difficult\'es du
th\'eor\`eme de Wintenberger.

Pour traiter la mauvaise r\'eduction justement, il me para\^{\i}t raisonnable de
postuler (\`a la suite de Kato) une th\'eorie de Dieudonn\'e log-cristalline
pour les log-sch\'emas en groupes finis localement libres.
L'anneau $\widehat{\Bst}$, variante de l'anneau
correspondant de Fontaine, mais poss\'edant une interpr\'etation en 
cohomologie log-cristalline devrait naturellement intervenir.
Une cons\'equence de ce travail en serait la semi-stabilit\'e 
(d\'ej\`a connue) des repr\'esentations galoisiennes
associ\'ees aux vari\'et\'es ab\'eliennes 
(du moins pour celles qui se prolongent en des log-vari\'et\'es ab\'eliennes),
avec le petit bonus que constitue la structure modulo $p^n$~;
en plongeant l'anneau $\widehat{\Bst}$
dans un anneau du type $\BdR$,
on devrait aussi retrouver le r\'esultat de Wintenberger.

\providecommand{\bysame}{\leavevmode ---\ }

\end{document}